\input amssym
\input xypic
\overfullrule=0pt
\def\Q{{\Bbb Q}}

\def\ds{{\frak{ds}}}
\def\dm{{\frak{dm}}}

\def\e{{\bf e}}
\def\b{{\bf b}}
\def\a{{\bf a}}
\def\c{{\bf c}}
\def\kv{{\frak{krv}_2}}
\def\tder{{\frak{tder}}}
\def\sder{{\frak{sder}}}
\centerline{\bf Double shuffle and Kashiwara-Vergne Lie algebras}
\vskip .3cm
\centerline{Leila Schneps}
\vskip 1cm
\noindent 
{\narrower{\narrower{{\bf Abstract.} We prove that the double shuffle Lie 
algebra $\ds$, dual to the space of new formal multiple zeta values, injects into the 
Kashiwara-Vergne Lie algebra $\kv$ defined and studied by Alekseev-Torossian.  The proof 
is based on a reformulation of the definition of $\kv$, and uses a theorem of 
Ecalle on a property of elements of $\ds$.}\par}\par}
\vskip 1cm
\noindent {\bf \S 1. Definitions and main results}
\vskip .5cm
Let $\Q\langle x,y\rangle$ denote the ring of polynomials in non-commutative variables 
$x$ and $y$, and ${\rm Lie}[x,y]$ the Lie algebra of Lie polynomials inside it.  For each $n\ge 1$, 
let $\Q_n\langle x,y\rangle$ (resp. ${\rm Lie}_n[x,y]$) denote the subspace of homogeneous polynomials 
(resp. Lie polynomials) of degree $n$.  For $k\ge 1$, let $\Q_{n\ge k}\langle x,y\rangle$ (resp.
${\rm Lie}_{n\ge k}[x,y]$) denote the space of polynomials (resp. Lie polynomials)
all of whose monomials are of degree $\ge k$, i.e. the direct sum of the
$\Q_n\langle x,y\rangle$ (resp. ${\rm Lie}_n[x,y]$) for $n\ge k$.

The main theorem of this paper gives an injective map between two Lie algebras
studied in the literature concerning formal multiple zeta values: the 
{\it double shuffle Lie algebra} $\ds$, investigated in papers by Racinet and
Ecalle amongst others (the associated graded  of $\ds$ is also studied in
papers by Zagier, Kaneko and others), and the Kashiwara-Vergne Lie algebra
introduced in work of Alekseev and Torossian (cf. [AT]).  We begin by
recalling the definitions of these two Lie algebras.  As vector spaces,
both are subspaces of the free Lie algebra ${\rm Lie}[x,y]$.

For any non-trivial monomial $w$ and polynomial $f\in \Q\langle x,y\rangle$,
we use the notation $(f|w)$ for the coefficient of the monomial $w$ in the
polynomial $f$, and extend it by linearity to polynomials $w$ without constant term.  
Set $y_i=x^{i-1}y$ for all $i\ge 1$; then all words ending in $y$ can be 
written as words in the variables $y_i$. The stuffle product $st(u,v)\in \Q\langle x,y\rangle$
of two such words $u$ and $v$
is defined recursively by
$$st(1,u)=st(u,1)=u\ \ {\rm and}\ \ st(y_iu,y_jv)=y_i\,st(u,y_jv)+y_j\,st(y_iu,v)+ y_{i+j}\,st(u,v).$$
\vskip .3cm
\noindent {\bf Definition 1.1.}  The {\it double shuffle Lie algebra } $\ds$\footnote{$^*$}{The equivalence
of the present definition with the usual definition introduced in [R] is proven in [CS], Theorem 2, which proves
that if a polynomial $f\in {\rm Lie}_n[x,y]$ has the property of the present definition, then 
$f+{{(-1)^{n-1}}\over{n}} (f|x^{n-1}y)y^n$ satisfies the stuffle relations for all pairs of words $u,v$ ending
in $y$. Since the words ending in $x$ are not involved in this condition, this is equivalent to the assertion
that $\pi_y(f)+{{(-1)^{n-1}}\over{n}}(f|x^{n-1}y)y^n$ satisfies stuffle, where $\pi_y(f)$ denotes the projection
of $f$ onto just its words ending in $y$.  This is the standard form of the defining property of elements
of $\ds$.} is
the vector space of elements $f\in {\rm Lie}_{n\ge 3}[x,y]$ such that
$$\bigl(f\,\bigl|\,st(u,v)\bigr)=0$$
for all words $u,v\in \Q\langle x,y\rangle$ ending in $y$ but not both simultaneously
powers of $y$.

\vskip .3cm
It has been shown by Racinet [R] (see also a simplified version 
of Racinet's proof in the appendix of [F]) and Ecalle [E] that $\ds$ is actually closed, i.e.
a Lie algebra, under the Poisson bracket defined on ${\rm Lie}[x,y]$ by
$$\{f,g\}=[f,g]+D_f(g)-D_g(f),\eqno(1.1)$$
where for any $f\in {\rm Lie}[x,y]$, the associated derivation $D_f$ of
${\rm Lie}[x,y]$ is defined by $D_f(x)=0$, $D_f(y)=[y,f]$.  This Lie
bracket corresponds to identifying $f$ with $D_f$ and taking the natural Lie
bracket on derivations:
$$[D_f,D_g]=D_{\{f,g\}}.\eqno(1.2)$$
\vskip .3cm
Let us now recall the definition of $\kv$.  Following [AT], let $TR$ denote
the vector space quotient of $\Q\langle x,y\rangle$ by relations
$ab=ba$.  The image in $TR$ of a monomial $w$ is the equivalence class of
monomials obtained by cyclically permutating the letters of $w$.  The
trace map $\Q\langle x,y\rangle\rightarrow TR$ is denoted by $tr$.
\vskip .3cm
\noindent {\bf Definition 1.2.} For any pair of elements $F$, $G\in 
{\rm Lie}_n[x,y]$ with $n\ge 1$, let $D_{F,G}$ denote the derivation of ${\rm Lie}[x,y]$
defined by $x\mapsto [x,G]$ and $y\mapsto [y,F]$.  Such a derivation
is said to be {\it special} if $x+y\mapsto 0$, i.e. if $[x,G]+[y,F]=0$.
The underlying vector space of the Kashiwara-Vergne Lie algebra is spanned by
those of these special derivations $D_{F,G}$ that also satisfy the
property that writing $F=F_xx+F_yy$ and $G=G_xx+G_yy$, there exists a
constant $A$ such that 
$$tr(F_yy+G_xx)\equiv A\,tr\bigl((x+y)^n-x^n-y^n\bigr)\in TR.\eqno(1.3)$$
It is shown in [AT] that $\kv$ is a Lie algebra under the 
natural bracket on derivations.  The degree provides a grading on $\kv$, for which
$(\kv)_n$ is spanned by the $D_{F,G}$ with $F,G\in {\rm Lie}_n[x,y]$. 
\vskip .3cm
The first graded piece, $(\kv)_1$, is 1-dimensional, generated by $D_{y,x}$.  The second graded
piece $(\kv)_2=0$.  Now let $n\ge 3$. Note that for any $F\in {\rm Lie}_n[x,y]$, if there exists $G
\in {\rm Lie}_n[x,y]$ such that $[y,F]+[x,G]=0$, then $G$ is unique.  
Indeed, $G$ is defined up to a centralizer of $x$, but that can only be
$x$, which is of degree $1$.  One of the most useful results of this paper is 
the precise determination of the elements $F$ admitting such a $G$, together 
with an explicit formula for $G$ (theorem 2.1, see also (1.5)).
\vskip .3cm
Let $\partial_x$ denote the derivation of 
$\Q\langle x,y\rangle$ defined by $\partial_x(x)=1$, $\partial_x(y)=0$.  
Following Racinet [R], for any polynomial $h$ in $x$ and $y$, set
$$s(h)=\sum_{i\ge 0} {{(-1)^i}\over{i!}}\partial_x^i(h)yx^i.\eqno(1.4)$$
Racinet shows that if $f=f_xx+f_yy$ is an element of ${\rm Lie}[x,y]$,
or indeed any polynomial such that $\partial_x(f)=0$, then
$$f=s(f_y).\eqno(1.5)$$
The main result of this paper is the following.
\vskip .5cm
\noindent {\bf Theorem 1.1.} {\it Let $\tilde f(x,y)\in \ds$, and set
$f(x,y)=\tilde f(x,-y)$ and $F(x,y)=f(z,y)$ with $z=-x-y$.  Write
$F=F_xx+F_yy=xF^x+yF^y$ in $\Q\langle x,y\rangle$.  Set $G=s(F^x)$.
Then the map $\tilde f\mapsto D_{F,G}$ yields an injective map of Lie algebras
$$\ds\hookrightarrow \kv.$$}
\vskip .05cm
\noindent {\bf Remark.} The map defined in theorem 1.1. from $\ds$ to the 
space of derivations $D_{F,G}$ mapping $x\mapsto [x,G]$ and $y\mapsto [y,F]$ 
is injective. Indeed, because $D_{F,G}(y)=[y,F]$ is a Lie element in which no
word starts and ends with $x$, we can recover $F$ from $D_{F,G}(y)$ 
by applying proposition 2.2 (with $x$ and $y$ exchanged in the statement),
and then we recover $\tilde f$ by $F(x,y)=\tilde f(z,-y)$.  

Furthermore, this injection of vector spaces is in fact an injection of
Lie algebras, since $\ds$ is equipped with the Poisson bracket, which is
compatible with the natural bracket on derivations (cf. (1.2)).  

Thus, to prove theorem 1.1, it remains only to prove that the derivations $D_{F,G}$ 
arising from elements $\tilde f\in\ds$ actually lie in $\kv$, i.e. are 
special and satisfy the trace formula (1.3).
\vskip .2cm
One of the main ingredients in our proof of theorem 1.1 is a combinatorial
reformulation of the defining properties of $\kv$, given in theorem 1.2
below.  First we need some definitions.
\vskip .3cm
\noindent {\bf Definition 1.3.}  Let $w=x^{a_0}y\cdots yx^{a_r}$ be a 
monomial in $\Q\langle x,y\rangle$ of depth $r$ (i.e. containing $r$
$y$'s), with $a_i\ge 0$ for $0\le i\le r$.  Let $anti$ denote the 
palindrome or backwards-writing operator on monomials, and let $push$ denote the cyclic
permutation of $x$-powers operator on monomials, defined respectively by
$$anti(x^{a_0}y\cdots yx^{a_{r-1}}yx^{a_r})=x^{a_r}yx^{a_{r-1}}y\cdots yx^{a_0}\eqno(1.6)$$
$$push(x^{a_0}y\cdots yx^{a_{r-1}}yx^{a_r})=x^{a_r}yx^{a_0}\cdots yx^{a_{r-1}}.\eqno(1.7)$$
For any word $w$, we define the list $Push(w)$ to be the list of $(r+1)$ words
obtained from $w$ by iterating the $push$ operator.  Note that $Push(w)$
is a list, not a set; it may contain repeated words.  For example, if
$w=x^2yxy$, then 
$Push(w)=[x^2yxy,yx^2yx,xy^2x^2]$, and if $w=xyxyx$, then $Push(w)=
[xyxyx,xyxyx,xyxyx]$.  
\vskip .3cm
\noindent {\bf Definition 1.4.} We extend the $anti$ and $push$ operators to operators
on polynomials by linearity; it makes sense to apply these operators to a polynomial 
even if the monomials in the polynomial have different degrees and depths.  If $f$ is a
polynomial in $x$ and $y$ of homogeneous degree $n\ge 3$, we say that $f$ is 
\vskip .1cm
\noindent $\bullet$ 
{\it palindromic} if $f=(-1)^{n-1}\,anti(f)$, 
\vskip .1cm
\noindent $\bullet$ {\it antipalindromic} if $f=(-1)^{n}\,anti(f)$,  
\vskip .1cm
\noindent $\bullet$ {\it push-invariant} if $push(f)=f$,
\vskip .1cm
\noindent $\bullet$ {\it push-constant} if there exists a constant $A$ such 
that $\sum_{v\in Push(w)} (f|v)=A$ for all $w\ne y^n$, and $(f|y^n)=0$.
\vskip .4cm
The following statement contains our reformulation of the definition of
$\kv$ that appears in [AT].
\vskip .4cm
\noindent {\bf Theorem 1.2.} {\it Let $V_{kv}$ be the vector space spanned by 
all polynomials $F\in {\rm Lie}_n[x,y]$ for $n\ge 3$ such that, writing $F=F_xx+F_yy$, we have
\vskip .1cm
\noindent
i) $F_y$ is antipalindromic, or equivalently, $F$ is push-invariant;
\vskip .1cm
\noindent
ii) $F_y-F_x$ is push-constant.
\vskip .1cm
\noindent
For each such $F$, set $G=s(F^x)$.  Then the map $F\mapsto D_{F,G}$ 
extends to a vector space isomorphism
$$V_{kv}\buildrel\sim\over\rightarrow \kv.\eqno(1.8)$$}
\vskip .3cm
The main result of \S 2, theorem 2.1, is an enumeration of several 
conditions equivalent to the specialness property.  
Using this result, theorems 1.1 and 1.2 are proved in \S 3.  
The proof of theorem 1.1 is based on two previously 
known results for $\ds$, each implying one of the two properties of theorem 
1.2.  The first of these theorems, theorem 3.3, is a translation into the
standard terms of $x,y$ variables of a theorem due to J. Ecalle [E]. Because
this result is couched in Ecalle's own original language, we give not only the 
reference to the precise statement, but also an appendix giving the complete 
calculation-translation which brings it to the form of theorem 3.3.  The second, 
theorem 3.4, appeared as Theorem 1 of [CS], with a complete elementary proof which 
was also based on an idea of Ecalle.
\vskip .8cm
\noindent {\bf \S 2. Characterizing special derivations}
\vskip .5cm
The main theorem of this section characterizes special derivations $D_{F,G}$
of ${\rm Lie}[x,y]$.  From now on, if $f$ is an element of ${\rm Lie}[x,y]$,
we say that $f$ is {\it special} if setting $F=f(z,y)$ with $z=-x-y$,
there exists a $G\in {\rm Lie}[x,y]$ such that $D_{F,G}$ is special.  By 
additivity, we may restrict ourselves to homogeneous Lie elements.  
\vskip .2cm
\noindent {\bf Notation.} For any $f\in \Q\langle x,y\rangle$, we will use the 
notation
$$f=f_xx+f_yy=xf^x+yf^y.$$
Observe that since every Lie element is palindromic, if
$f\in {\rm Lie}_n[x,y]$, we have
$$f=(-1)^{n-1}\,anti(f)=f_xx+f_yy=(-1)^{n-1}x\,anti(f_x)+(-1)^{n-1}y\,anti(f_y)=
xf^x+yf^y,$$
so in fact
$$f^x=(-1)^{n-1}\,anti(f_x),\ \ \ \ f^y=(-1)^{n-1}\,anti(f_y).\eqno(2.1)$$
Recall also the definition of the map $s:\Q\langle x,y\rangle
\rightarrow \Q\langle x,y\rangle$ from (1.4).  We will also use the similar 
map
$$s'(h)=\sum_{i\ge 0} {{(-1)^i}\over{i!}} x^iy\partial_x^i(h).\eqno(2.2)$$
When $f\in {\rm Lie}_n[x,y]$ ($n\ge 2$), it follows by symmetry from Racinet's result $f=s(f_y)$ 
that if we write $f=xf^x+yf^y$, then $f=s'(f^y)$.
\vskip .3cm
\noindent {\bf Theorem 2.1.} {\it Let $n\ge 3$, and let $f\in {\rm Lie}_n[x,y]$.
Set $F=f(-x-y,y)$, and write $f=f_xx+f_yy$ and $F=F_xx+F_yy$.
Then the following are equivalent:
\vskip .1cm
i) $f$ is special, i.e. there exists a unique $G\in {\rm Lie}_n[x,y]$ such
that $[y,F]+[x,G]=0$.
\vskip .1cm
ii) Setting $G=s'(F_x)$, the derivation $D_{F,G}$ is special.
\vskip .1cm
iii) $F_y$ is antipalindromic.
\vskip .1cm
iv) $F$ is push-invariant.
\vskip .1cm
v) $f_y-f_x$ is antipalindromic.}
\vskip .4cm
The equivalence of i), ii) and iii) is given in proposition 2.3.
the equivalence of iii) and iv) is proven in proposition 2.4, and the 
equivalence of iii) and v) is given following proposition 2.6. Some of these
results, in particular propositions 2.2 and 2.6, will also be used in the proofs
of the main theorems in \S 3.  

\vskip .3cm
\noindent {\bf Proposition 2.2.} {\it Let $n\ge 3$, and let
$f\in {\rm Lie}_n[x,y]$ have the property that expanded as a polynomial, 
$f$ has no terms that start and end in $y$, so that writing $f=f_xx+f_yy$, 
we have $f_yy=xPy$. Then $s(P)\in {\rm Lie}_{n-1}[x,y]$ and $f=[x,s(P)]$.}
\vskip .2cm
\noindent {\bf Proof.} By hypothesis, $f$ has no terms starting and ending in $y$,
so we can write $f_yy=xPy$.
By Racinet's result, we have $g=s(g_y)$ for all $g\in {\rm Lie}_n[x,y]$ with $n\ge 2$, 
so in particular we have $f=s(xP)$.  Now, since the partial derivative
satisfies $\partial^i(xP)=i\partial^{i-1}(P)+x\partial^i(P)$,
and $\partial^n(P)=0$ since $P$ is of degree $n-1$, we compute
$$\eqalign{f=s(xP)&=\sum_{i=0}^n {{(-1)^i}\over{i!}}\partial^i(xP)yx^i\cr
&=\sum_{i=0}^n {{(-1)^i}\over{i!}}(i\partial^{i-1}(P)+x\partial^i(P))yx^i\cr
&=\sum_{i=0}^n {{(-1)^i}\over{(i-1)!}}\partial^{i-1}(P)yx^i
+\sum_{i=0}^n {{(-1)^i}\over{i!}}(x\partial^i(P))yx^{i+1}\cr
&=\sum_{i=0}^{n-1} {{(-1)^{i-1}}\over{i!}}\partial^i(P)yx^i
+\sum_{i=0}^{n-1} {{(-1)^i}\over{i!}}(x\partial^i(P))yx^i\cr
&=-s(P)x+x\,s(P).}$$
Thus, $f=[x,s(P)]$.

It remains only to show that $s(P)$ is a Lie element. 
Let $\Phi:\Q_{n\ge 1}\langle x,y
\rangle\rightarrow {\rm Lie}[x,y]$ be the linear map sending a non-trivial word
$w=x_1x_2x_3\cdots x_m$ to $[x_1,[x_2,[x_3,\cdots]]]$, where $x_i\in \{x,y\}$,
and let $\theta:\Q\langle x,y\rangle\rightarrow {\rm End}_\Q{\rm Lie}[x,y]$
be the algebra homomorphism mapping $x$ to $ad(x)$ and $y$ to $ad(y)$.
By [B, Ch 2, \S 3, no. 2] the following properties hold:
\vskip .1cm
$\bullet$ a polynomial $h\in \Q_n\langle x,y\rangle$ is Lie if and only if
$\Phi(h)=nh$;
\vskip .1cm
$\bullet$ $\Phi(uv)=\theta(u)\Phi(v)$ for $u\in \Q\langle x,y\rangle$
and $v\in \Q_{n\ge 1}\langle x,y\rangle$.
\vskip .1cm
$\bullet$ $\theta(u)(v)=[u,v]$ if $u$ is Lie.
\vskip .1cm
Since $f\in {\rm Lie}[x,y]$, we have
$$[f,x]=\theta(f)(x)=\theta([x,s(P)])(x)=\bigl[ad(x),\theta\bigl(s(P)\bigr)\bigr](x)=
\bigl[x,\theta\bigl(s(P)\bigr)(x)\bigr]=-\bigl[\theta\bigl(s(P)\bigr)(x),x\bigr].$$
Thus, $\bigl[f+\theta\bigl(s(P)\bigr)(x),x\bigr]=0$, so since both $f$ and 
$\theta\bigl(s(P)\bigr)(x)$ are Lie elements of degree $>1$, we have 
$f=-\theta\bigl(s(P)\bigr)(x)$.  Thus, 
$$nf=\Phi(f)=\Phi([x,s(P)])=\theta(x)\Phi\bigl(s(P)\bigr)-\theta\bigl(s(P)\bigr)\Phi(x)=
\bigl[x,\Phi\bigl(s(P)\bigr)\bigr]-\theta\bigl(s(P)\bigr)(x)=\bigl[x,\Phi\bigl(s(P)\bigr)\bigr]+f.$$
Thus $\bigl[x,\Phi\bigl(s(P)\bigr)\bigr]=(n-1)f=(n-1)[x,s(P)]$, so 
$\bigl[x,\Phi\bigl(s(P)\bigr)-(n-1)s(P)\bigr]=0$.  Since $s(P)$ is of degree $n-1>1$, we must have
$\Phi\bigl(s(P)\bigr)=(n-1)s(P)$, but this means that $s(P)\in {\rm Lie}_{n-1}[x,y]$.
\hfill{$\diamondsuit$}
\vskip .3cm
\noindent {\bf Proposition 2.3.} {\it Let $f\in {\rm Lie}_{n\ge 3}[x,y]$,
and set $F=f(z,y)=F_xx+F_yy$ and $G=s'(F_x)$.  Then $D_{F,G}$ is special
if and only if $f$ is special, and this is the case if and only if
$F_y$ is antipalindromic.}
\vskip .2cm
\noindent {\bf Proof.} If setting $G=s'(F_x)$, the derivation $D_{F,G}$ is
special, then $f$ is special by definition.  Conversely, if $f$ is special,
there exists a unique $G\in {\rm Lie}_n[x,y]$ such that $[y,F]+[x,G]=0$. 
Setting $H=yF-Fy=Gx-xG$ and writing 
$F=F_xx+F_yy=xF^x+yF^y$ and $G=G_xx+G_yy=xG^x+yG^y$, this means that
$$H=yF_yy+yF_xx-yF^yy-xF^xy=xG^xx+yG^yx-xG_xx-xG_yy,\eqno(2.3)$$
so comparing the terms starting with $x$ and ending with $y$, we find that
$-xF^xy=-xG_yy$, so $F^x=G_y$.  By a result of Racinet [R], since $G$ is
a Lie element, we must have $G=s(G_y)=s(F^x)=s'(F_x)$.  This proves the
first equivalence.

Let us now assume that $F_y$ is antipalindromic, i.e. by (2.1), $F^y=F_y$. 
Set 
$$H=yF-Fy=y(F_yy+F_xx)-(yF^y+xF^x)y=yF_yy-yF^yy+yF_xx-xF^xy.\eqno(2.4)$$
This shows that $H$ has no words starting and ending in $y$, so
by proposition 2.2, there exists $G\in {\rm Lie}_{n-1}[x,y]$ 
such that $H=Gx-xG$.  But then the derivation $D_{F,G}$ is special,
so $f$ is special. 

Finally, assume that $f$ is special, and set $H=yF-Fy$, so that there exists
$G$ with $H=yF-Fy=Gx-xG$. Then (2.3) holds.  The expression $H=Gx-xG$ shows that
$H$ can have no terms starting and ending in $y$, and the left-hand expression
for $H$ in (2.3) then shows that we must have $F_y=F^y$,
i.e. by (2.1), $F_y$ is antipalindromic.\hfill{$\diamondsuit$}
\vskip .4cm
\noindent {\bf Proposition 2.4.} {\it Let $F\in {\rm Lie}_n\langle x,y\rangle$.
Then $F_y$ is antipalindromic if and only if $F$ is push-invariant.}
\vskip .2cm
\noindent {\bf Proof.} As usual, we write $F=F_xx+F_yy=xF^x+yF^y$.  Assume
first that $F_y$ is antipalindromic, i.e. that
$F_y=F^y$.  Since $F$ is a Lie
polynomial, we have $F=s(F_y)=s'(F^y)=s'(F_y)$, i.e.
$$F=\sum_{i\ge 0} {{(-1)^i}\over{i!}} \partial_x^i(F_y)yx^i=
\sum_{i\ge 0} {{(-1)^{i}}\over{i!}} x^iy\partial_x^i(F^y)
=\sum_{i\ge 0} {{(-1)^{i}}\over{i!}} x^iy\partial_x^i(F_y).\eqno(2.5)$$ 
Using the second and fourth terms of (2.5), we compute
the coefficient of a word in $F$ as
$$\eqalign{(F|x^{a_0}y\cdots x^{a_{r-1}}yx^{a_r})&=
{{(-1)^{a_r}}\over{(a_r)!}}(\partial_x^{a_r}(F_y)yx^{a_r}|
x^{a_0}y\cdots yx^{a_{r-1}}yx^{a_r})\cr
&={{(-1)^{a_r}}\over{(a_r)!}}(\partial_x^{a_r}(F_y)|
x^{a_0}y\cdots yx^{a_{r-1}})\cr
&={{(-1)^{a_r}}\over{(a_r)!}}(x^{a_r}y\partial_x^{a_r}(F_y)|
x^{a_r}yx^{a_0}y\cdots yx^{a_{r-1}})\cr
&=(F|x^{a_r}yx^{a_0}y\cdots yx^{a_{r-1}}),}$$
so $F$ is push-invariant.

In the other direction, suppose that $F$ is push-invariant, and let's show that 
$F_y=F^y$.  By assumption, we have
$$(F|x^{a_0}y\cdots yx^{a_r})=(F|x^{a_r}yx^{a_0}y\cdots x^{a_{r-1}}).$$
In particular, for all words with $a_r=0$, we have
$(F|x^{a_0}y\cdots yx^{a_{r-1}}y)=(F|yx^{a_0}y\cdots yx^{a_{r-1}}),$
i.e.
$$(F_yy|x^{a_0}y\cdots yx^{a_{r-1}}y)=(yF^y|yx^{a_0}y\cdots yx^{a_{r-1}}),$$
so
$$(F_y|x^{a_0}y\cdots yx^{a_{r-1}})=(F^y|x^{a_0}y\cdots yx^{a_{r-1}}).$$
Thus $F_y=F^y$.  \hfill{$\diamondsuit$}
\vskip .3cm
\noindent {\bf Lemma 2.5.} {\it Let $g\in \Q_n\langle x,y\rangle$, let
$\phi(x,y)$ and $\psi(x,y)$ be linear expressions of the form $ax+by$, $a,b\in \Q$,
and let $h(x,y)=g\bigl(\phi(x,y),\psi(x,y)\bigr)$.  If $g$ is antipalindromic, 
then $h$ is antipalindromic.}
\vskip .2cm
\noindent {\bf Proof.} The operator $anti$ is an anti-automorphism of the
ring $\Q\langle x,y\rangle$, so 
$$anti(h)=anti\Bigl(g\bigl(anti(\phi),anti(\psi)\bigr)\Bigr).$$
But $anti$ fixes linear expressions $ax+by$, so since $g$ is antipalindromic, we have
$$anti(h)=anti\Bigl(g\bigl(\phi,\psi\bigr)\Bigr)=anti(g)(\phi,\psi)
=(-1)^{n-1}g(\phi,\psi)=(-1)^{n-1}h.$$
Thus $h$ is antipalindromic.\hfill{$\diamondsuit$}
\vskip .3cm
\noindent {\bf Proposition 2.6.} {\it For any $g\in {\rm Lie}_n[x,y]$,
set $z=-x-y$ and $G=g(z,y)$. Write $g=g_xx+g_yy$ and $G=G_xx+G_yy$.  Then
$$G_y-G_x=g_y(z,y).$$
In particular, $g_y$ is antipalindromic if and only if $G_y-G_x$ is 
antipalindromic.}
\vskip .2cm
\noindent {\bf Proof.} We have $g(x,y)=g_x(x,y)x+g_y(x,y)y$, so
$$G=g(z,y)=g_x(z,y)z+g_y(z,y)y=-g_x(z,y)x-g_x(z,y)y+g_y(z,y)y.$$
Thus $G_y=-g_x(z,y)+g_y(z,y)$ and $G_x=-g_x(z,y)$, so $G_y-G_x=g_y(z,y)$.
Then by Lemma 2.5, since $g_y$ is antipalindromic, so is $G_y-G_x$, and
the converse holds as well since $(G_y-G_x)(z,y)=g_y$.\hfill{$\diamondsuit$} 
\vskip .5cm
We can now conclude the proof of theorem 2.1 by showing the equivalence of iii) and v).
To do this, we simply apply proposition 2.6 with $f=G$ and $g=F$, to see that
$F_y$ is antipalindromic if and only if $f_y-f_x$ is antipalindromic.
This completes the proof. 
\vskip .8cm
\noindent {\bf \S 3. Proofs of theorems 1.2 and 1.1}
\vskip .5cm
\noindent {\bf Proof of theorem 1.2.}  Let $F\in V_{kv}$.  We may assume
that $F$ is homogeneous of degree $n\ge 3$, i.e. $F\in {\rm Lie}_n[x,y]$ with
$n\ge 3$.
Set $G=s'(F_x)=s(F^x)$.  By theorem 2.1, $F_y$ is antipalindromic if
and only if $F$ is push-invariant, and these conditions are equivalent
to the fact that $G\in {\rm Lie}_n[x,y]$ and $D_{F,G}$ is special.

Now consider the map $F\mapsto D_{F,G}$ from $V_{kv}$ to the vector
space of special derivations, and let us show that it is injective.  Suppose that
$F,F'\in V_{kv}$ and $D_{F,G}=D_{F',G'}$.  Then $D_{F,G}(y)=D_{F',G'}(y)$,
i.e. $[y,F]=[y,F']$, so $F-F'$ commutes with $y$.  Since $F-F'$ is of degree $>1$,
this means that $F-F'=0$.

Let us now show that $D_{F,G}$ satisfies the trace formula (1.3).
Note that by (2.1), $F_y-F_x=(-1)^{n-1}anti(F^y-F^x)$, so by
symmetry, $F_y-F_x$ is push-constant if and only if $F^y-F^x$ is 
push-constant.  It is convenient to use the latter condition. 

Since any Lie polynomial of degree $>1$ is a sum of terms of the form
$fg-gf$, Lie polynomials map to zero in TR.  Thus, we have 
$tr(G_xx)=-tr(G_yy)$, so 
$$\eqalign{tr(F_yy+G_xx)&=tr(F_yy-G_yy)\cr
&=tr(F_yy-F^xy)\ \ \ {\rm since}\ G_y=F^x\cr
&=tr(F^yy-F^xy)\ \ \ {\rm since}\ F^y=F_y\cr
&=tr\bigl((F^y-F^x)y\bigr).}\eqno(3.1)$$
Rephrasing the trace formula (1.3) via (3.1) as 
$$tr((F^y-F^x)y)=A\,tr\bigl((x+y)^n-x^n- y^n\bigr),\eqno(3.2)$$ 
we can now show that a special derivation $D_{F,G}$ satisfies the trace formula
in $TR$ if and only if $F^y-F^x$ is push-constant.  In fact,
these are just two ways of making the identical statement.  
To see this, let $\overline C$ denote the {\it list} of words in the cyclic permutation class of $w$,
so that $\overline C$ contains exactly $n$ words; then $\overline C$ consists of $n/|C|$
copies of $C$.  For any word $v=uy$ ending in $y$, let $\overline C$ denote the associated cyclic 
permutation list, and $\overline C_y$ the list obtained from $\overline C$ by removing all the words 
ending in $x$.  Write $\overline C_y=[u_1y,\ldots,u_ry]$.  Then by definition, we have the 
equality of lists
$$[u_1,\ldots,u_r]=Push(u).\eqno(3.3)$$
Now, the trace condition $tr((F^y-F^x)y)=A\,tr\bigl((x+y)^n-x^n-y^n\bigr)$ means firstly
that $\bigl((F^y-F^x)y|y^n\bigr)=0$, which is equivalent to $\bigl(F^y-F_x|y^{n-1})=0$,
and secondly that for each equivalence class $C$ of 
cyclic permutations of a given word $w\ne y^n$, we have
$$\Bigl(tr\bigl((F^y-F^x)y\bigr)\,\bigl|\,C\Bigr):=\sum_{v\in C} \Bigl((F^y-F^x)y\,\bigl|\,v\Bigr)=
|C|A,$$
where the first equality is just the definition of the coefficient of an equivalence class in a trace 
polynomial. Using the notation $\overline C$ and $\overline C_y$ as above and (3.3), this means that
for every word $v\ne y^n$ ending in $y$, writing $v=uy$, we have 
$$|C|A={{|C|}\over{n}}\sum_{v\in \overline C} \Bigl((F^y-F^x)y\,\bigl|\,v\Bigr)
={{|C|}\over{n}}\sum_{v\in \overline C_y} \bigl((F^y-F^x)y\,\bigl|\,v\bigr)
={{|C|}\over{n}}\sum_{u'\in Push(u)} \bigl((F^y-F^x)\,\bigl|\,u'\bigr).
\eqno(3.4)$$

But this is equivalent to
$$\sum_{u'\in Push(u)} \bigl((F^y-F^x)\,\bigl|\,u'\bigr)=nA\eqno(3.5)$$
for all $u'\ne y^{n-1}$,
which, together with the fact that $(F^y-F^x|y^{n-1})=0$, is precisely 
equivalent to the statement that $F^y-F^x$ is push-constant (for the constant $nA$).

So far we have proven that $F\mapsto D_{F,s(F^x)}$ for homogeneous $F$ extends to 
an injective map $V_{kv}\hookrightarrow \kv$. Let us show that it is an isomorphism, i.e.
also surjective.  It is enough to consider derivations $D_{F,G}\in \kv$ with $F$, $G$ homogeneous 
of degree $n$.  Then $D_{F,G}$ is special, so $F_y$ is antipalindromic by theorem 2.1, and $D_{F,G}$ 
satisfies the trace formula (3.2), which as we just saw is equivalent to the property that 
$F^y-F^x$ is push-constant.  Finally, since $F_y-F_x=(-1)^{n-1}anti(F^y-F^x)$,
we see that $F_y-F_x$ is also push-constant, so $F\in V_{kv}$, completing
the proof.\hfill{$\diamondsuit$} 
\vskip .4cm
Let us now prove theorem 1.1.  The proof is based on the fact that two 
previously known combinatorial results about double shuffle elements 
$\tilde f\in \ds$ make it possible to deduce that $F=\tilde f(x,-y)$ satisfies
the two defining properties of $V_{kv}$ given in theorem 1.2.  Thus 
$\tilde f\mapsto F$ yields an injection $\ds\hookrightarrow
V_{kv}$, and the injection $V_{kv}\hookrightarrow \kv$ of theorem 1.2
completes the argument.  

The two known results are
given in theorems 3.3 and 3.4.  As the original statement of theorem 3.3
is extremely different in appearance (theorem A.1 below), the translation from 
the original terminology to the statement given here is provided in the 
appendix, which also serves as an initiation to Ecalle's language.   We
write $\ds_n$ for the homogeneous weight $n$ part of $\ds$, consisting of
polynomials in $\ds$ which are of homogeneous degree $n$.
\vskip .3cm
\noindent {\bf Theorem 3.3.} [E, cf. Appendix] {\it Let $\tilde f\in\ds_n$,
and write $\tilde f=\tilde f_xx+
\tilde f_yy$.  Then $\tilde f_x+\tilde f_y$ is antipalindromic.}
\vskip .3cm
\noindent {\bf Theorem 3.4.} [CS] {\it Let $\tilde f=\tilde f_xx+
\tilde f_yy\in \ds_n$, and set $A=(\tilde f|x^{n-1}y)$.  Then $\tilde f_y$ 
satisfies the property that $(\tilde f_y|y^{n-1})=0$ and for
each degree $n$ monomial $w\ne y^{n-1}$ containing $r$ $y$'s, we have
$$\sum_{v\in Push(w)} \bigl(\tilde f_y|v\bigr)=(-1)^rA.$$}
\vskip .3cm
\noindent {\bf Proof of theorem 1.1.} 
Let $n\ge 3$ and assume that $\tilde f\in \ds_n$, i.e. $\tilde f$ is a 
homogeneous Lie polynomial of degree $n$.  Set $f(x,y)=\tilde f(x,-y)$.  
It follows directly from theorem 3.4 that $f_y$ is push-constant. Let us
deduce from theorem 3.3 that $f_y-f_x$ is antipalindromic. Indeed,
$f(x,y)=f_x(x,y)x+f_y(x,y)y$ and $\tilde f(x,y)=f(x,-y)$, so
$\tilde f(x,y)=f_x(x,-y)x-f_y(x,-y)y$, i.e.
$\tilde f_x=f_x(x,-y)$, $\tilde f_y=-f_y(x,-y)$. Thus
$$\tilde f_x+\tilde f_y=f_x(x,-y)-f_y(x,-y)=(f_x-f_y)(x,-y).\eqno(3.6)$$
The left-hand side is antipalindromic by theorem 3.3, so the right-hand
side is antipalindromic, and then by Lemma 2.5 $f_x-f_y$ and thus
also $f_y-f_x$ are antipalindromic.

Set $F=f(z,y)$.  We will use the two properties on $f$ to show that
$f\mapsto F$ is an injection from $\ds$ into $V_{kv}$.  By proposition 2.6
with $g=F$ and $G=f$, we see that $f_y-f_x$ antipalindromic implies that
$F_y$ is antipalindromic.  It remains only to show that 
$f_y$ push-constant implies that $F_y-F_x$ is push-constant,
which is a little more delicate. We prove it in the following lemma.
\vskip .3cm
\noindent {\bf Lemma 3.5.} {\it For any $f\in {\rm Lie}_n[x,y]$,
set $F=f(z,y)$ and write $f=f_xx+f_yy$ and
$F=F_xx+F_yy$.  Suppose that $f_y$ is push-constant for a constant 
$A$, and that $A=0$ if $n$ is even.  Then $F_y-F_x$ is also push-constant for 
$A$.}
\vskip .2cm
\noindent {\bf Proof.} 
To show that $F_y-F_x$ is push-constant for $A$, let us first show
that $(F_y-F_x|y^{n-1})=0$.  As we saw in the proof of theorem 1.2, the 
condition that $f_y$ is push-constant is equivalent to the condition that 
$tr(f_yy)=A\,tr\bigl((x+y)^n-x^n-y^n\bigr)$.  By proposition 2.6 with
$g=f$ and $G=F$, we have
$F_y-F_x=f_y(z,y)$, so
$(F_y-F_x)(z,y)=f_y$.  Multiplying by $y$ on the right of both sides
and taking the trace yields
$$tr\Bigl((F_y-F_x)(z,y)y\Bigr)=tr(f_yy)
=A\,tr\Bigl((x+y)^n-x^n-y^n\Bigr).$$
Making the variable change $x\mapsto z=-x-y$ on both sides, this gives
$$tr\Bigl((F_y-F_x)y\Bigr)
=A\,tr\Bigl((-1)^nx^n-(-1)^n(x+y)^n-y^n\Bigr).$$
When $n$ is odd, the right-hand side does not contain the equivalence class
of $y^n$, so the left-hand side cannot contain it either, which means that
$(F_y-F_x\,|\,y^{n-1})=0$.  When $n$ is even, $A=0$ by assumption,
so the equivalence class of $y^{n-1}$ cannot appear in the left-hand side,
which again means that $(F_y-F_x\,|\,y^{n-1})=0$.

Now let us prove that $F_y-F_x$ is push-constant.  Write
$$f_y=\sum_\a c_\a\,x^{a_0}y\cdots yx^{a_r}=
\sum_v c_v\,v,$$
where $\a$ runs over the tuples $\a=(a_0,\ldots,a_r)$ with 
$r\ge 1$ and $a_0+\cdots+a_r=n-r-1$, and $v$ runs over degree $n-1$ words.
If $v=x^{a_0}y\cdots yx^{a_r}$, we write $c_v=c_\a$.  For a given
tuple $\a=(a_0,\ldots,a_r)$, let 
$$Push(\a)=[(a_0,\ldots,a_r),(a_r,a_0,\ldots,a_{r-1}),\ldots,
(a_1,\ldots,a_r,a_0)]$$
be the list of its $r+1$ cyclic permutations.
The fact that $f_y$ is push-constant means that for all $w\ne y^{n-1}$, we have
$$\sum_{v\in Push(w)} (f_y|v)=
\sum_{v\in Push(w)} c_v = \sum_{\a'\in Push(\a)} c_{\a'} = A.\eqno(3.7)$$

Let us now compute the coefficient in $F_y-F_x$ of a given word 
$w=x^{b_0}y\cdots yx^{b_d}$, $w\ne y^{n-1}$.  
By proposition 2.6, we have
$$\eqalign{F_y-F_x=f_y(z,y)&=\sum_\a c_\a\,z^{a_0}y\cdots yz^{a_r}
= \sum_\a (-1)^{n-r-1}c_\a\,(x+y)^{a_0}y\cdots y(x+y)^{a_r},}\eqno(3.8)$$
so
$$\bigl(F_y-F_x\,\bigl|\,w\bigr)= (-1)^{n-1}\biggl(\sum_\a 
(-1)^rc_\a\,(x+y)^{a_0}y\cdots y(x+y)^{a_r}\Bigl|\, w\biggr).
\eqno(3.9)$$
Clearly if $r>d$ then the expansion of $(x+y)^{a_0}y\cdots y(x+y)^{a_r}$ cannot
contain the word $w$, so (3.9) is equal to
$$\bigl(F_y-F_x\,\bigl|\,w\bigr)=(-1)^{n-1}\sum_{{\bf a}\ s.t.\ 0\le r\le d} 
(-1)^rc_{\bf a}\, \biggl((x+y)^{a_0}y\cdots y(x+y)^{a_r} \,
\Bigl|\, x^{b_0}y\cdots yx^{b_d}\biggr).\eqno(3.10)$$
The only terms $(x+y)^{a_0}y\cdots (x+y)^{a_r}$ in which $w$ will appear
with a positive coefficient (necessarily equal to $1$) are the
$2^d$ terms $(x+y)^{a_0}y\cdots (x+y)^{a_r}$ constructed as follows:
choose any of the $2^d$ subsets of the $y$'s in $w$,
and change the $y$'s in that subset to $x$'s; then substitute 
$x\mapsto (x+y)$ in the resulting word.

Let $w=x^{b_0}y\cdots yx^{b_d}$ be a monomial, and set $\b=(b_0,\ldots,b_d)$.  Write $X_{\bf b}$ 
for the set of $2^d$ sequences $(a_0,\ldots,a_r)$, 
$0\le r\le d$, such that the corresponding word $x^{a_0}y\cdots yx^{a_r}$ is 
obtained from $w$ by 
changing any subset of $y$'s into $x$'s.  Then the coefficient (3.10) is equal 
to
$$\bigl(F_y-F_x\,\bigl|\,w\bigr)=(-1)^{n-1}\sum_{{\bf a}\in X_{\bf b}} (-1)^r
c_{\bf a}.  \eqno(3.11)$$
By (3.11), we have
$$\eqalign{\sum_{v\in Push(w)} \bigl(F_y-F_x|v\bigr)&=
(-1)^{n-1}\sum_{{\bf c}\in Push(\b)}\sum_{{\bf a}\in X_{\bf c}} (-1)^rc_\a}.\eqno(3.12)$$
Let us write
$${\cal X}_\b=\coprod_{\c\in Push(\b)} X_\c$$
for the disjoint union, i.e. the list-union of the words in the lists $X_\c$,
where $\c$ runs through the cyclic permutations of $\b$.  There are $(d+1)2^d$
words in ${\cal X}_\b$.  Let us count the words in ${\cal X}_\b$ of each 
given depth $0\le r\le d$.

For each tuple $\c\in Push(\b)$, let $w_\c$ be the word associated to $\c$.
The list ${\cal X}_\b$ is exactly the list of all words obtained by changing 
$k$ of the $d$ $y$'s in $w_\c$ to $x$'s, for all $0\le k\le d$ and
all $\c\in Push(\b)$.  
Thus, ${\cal X}_\b$ contains $(d+1)$ words of depth $d$, which are the words 
$w_\c$ for $\c\in Push(\b)$, and for each smaller depth $r=d-k$ for 
$1\le k\le d$, ${\cal X}_\b$ contains the words obtained by changing 
$k$ $y$'s to $x$'s in each of the $d+1$ words (all of depth $d$) of $Push(\b)$. Thus, there
are exactly $(d+1)\bigl({{d}\atop{k}}\bigr)$ words of depth $r=d-k$
in ${\cal X}_\b$, and these words fall into exactly
$${{d+1}\over{d-k+1}}\Bigl({{d}\atop{k}}\Bigr)=\Bigl({{d+1}\atop{k}}\Bigr)$$ 
cycles of length $r+1=d-k+1$, of words of depth $r=d-k$.

Since $f_y$ is push-constant, the coefficients $c_\a$ of each of the 
$\bigl({{d+1}\atop{k}}\bigr)$ cycles of depth $r=d-k$ in 
$f_y$ add up to $A$.  Thus, for all $\b\ne (1,\ldots,1)$, (3.12) is given by
$$(-1)^{n-1}\sum_{\a\in {\cal X}_b} (-1)^rc_\a=(-1)^{n-1}\sum_{k=0}^{d} \Bigl({{d+1}\atop{k}}\Bigr)(-1)^{d-k}A=(-1)^n
\sum_{k=0}^{d} \Bigl({{d+1}\atop{k}}\Bigr)(-1)^{d+1-k}A=(-1)^{n-1}A.$$
This proves that $F_y-F_x$ is push-constant for the value $(-1)^{n-1}A$. 
\hfill{$\diamondsuit$}
\vskip .3cm
We can now conclude the proof of theorem 1.1. Using the well-known result on 
$\ds$ (cf. [E], [R], [IKZ]...) that the coefficient of $x^{n-1}y$ is zero for 
all even-degree elements of $\ds$, we see that when $n$ is even, 
$A=0$ in theorem 3.4, so if $\tilde f\in \ds$, then $f=\tilde f(x,-y)$ 
satisfies the hypotheses of Lemma 3.5.  

Thus, we have shown so far that if $\tilde f\in \ds$, setting 
$f(x,y)=\tilde f(x,-y)$ and $F=f(z,y)$, $F_y$ is antipalindromic by the argument
of of the first paragraph of the proof of theorem 1.1, 
and $F_y-F_x$ is push-constant by Lemma 3.5.
Thus, the map $\tilde f\mapsto F$ is an injective map from $\ds\rightarrow
V_{kv}$.  By (1.8), we then have an injective composition of maps
$$\ds\hookrightarrow V_{kv}\buildrel\sim\over\rightarrow \kv.$$
This concludes the proof of theorem 1.1.\hfill{$\diamondsuit$}
\vfill\eject
\noindent {\bf \S 4. The prounipotent version}
\vskip .5cm
Let $V$ be a graded vector space, and let $\frak{un}(V)$ denote the Lie
algebra of (pro)unipotent endomorphisms of $V$, i.e. linear endomorphisms
$D$ such that $D(V_{\ge n})\subset V_{\ge n+1}$.  The usual exponentiation
$$exp(D)=\sum_{n\ge 0} {{1}\over{n!}}D^n\eqno(4.1)$$
maps $\frak{un}(V)$ bijectively to the group $UN(V)$ of (pro)unipotent linear
automorphisms of $V$.

Suppose we now have a Lie algebra $\frak{g}$ equipped with an injective
Lie algebra map $\frak{g}\buildrel\rho\over\hookrightarrow \frak{un}(V)$.  The universal enveloping 
algebra ${\cal U}\frak{g}$ is a ring whose multiplication we denote by
$\odot$. The exponential associated to $\frak{g}$ is given by the formula
$$exp^\odot(f)=\sum_{n\ge 0} {{1}\over{n!}}f^{\odot n};\eqno(4.2)$$
it maps $\frak{g}$ bijectively to the associated group $G\subset
\widehat{{\cal U}\frak{g}}$, and the following diagram commutes:
$$\xymatrix{G\ar[r]&UN(V)\\
\frak{g}\ar[r]^\rho\ar[u]^{{\rm exp}^\odot}&\frak{un}(V).\ar[u]_{\rm exp}}\eqno(4.3)$$

\vskip .2cm
Now let $V$ denote the underlying vector space of ${\rm Lie}[x,y]$. Following 
the notation of [AT], let $\tder_2$ denote the Lie algebra of {\it tangential
derivations} of ${\rm Lie}[x,y]$, i.e. derivations $D$ having the property 
that $D(x)=[x,a]$ and $D(y)=[y,b]$ for elements $a,b\in {\rm Lie}[x,y]$. 
There is an injective map of Lie algebras $\tder_2\hookrightarrow \frak{un}(V)$.
Indeed, if $V={\rm Lie}[x,y]$ is equipped with the grading given by the 
degree, then any derivation $D\in \tder_2$ increases the degree, i.e. 
$D(V_{\ge n})\subset V_{\ge n+1}$.
Let $TAut_2$ denote the group of automorphisms of $V$ obtained
by exponentiating $\tder_2$:
$$\eqalign{{\rm exp}:\tder_2&\rightarrow TAut_2\subset UN(V)\cr
D&\mapsto {\rm exp}(D)=\sum_{n\ge 0} {{1}\over{n!}} D^n.}\eqno(4.4)$$
%If $D=D_f$, then by (4.4), setting $F={\rm exp}^\odot(f)$, we have
%$$A(x)=F^{-1}(z,x)\,x\,F(z,x),\ \ \  
%A(y)=F^{-1}(z,y)\,y\,F(z,y).\eqno(4.5)$$
\vskip .1cm
Let $\sder_2$ denote the subalgebra
of $\tder_2$ consisting of derivations $D$ such that $D(x+y)=0$, and $SAut_2$ 
the corresponding subgroup of $TAut_2$ consisting of automorphisms such that 
$A(x+y)=x+y$, so that ${\rm exp}(\sder_2)=SAut_2$.  According to [AT], the 
exponential map (4.1) not only restricts to (4.4), but also
to maps from the following subspaces to subgroups:
$$\xymatrix{KRV_2\ar[r]&SAut_2\ar[r]&TAut_2\ar[r]&UN(V)\\
\kv\ar[r]\ar[u]^{\rm exp}&\sder_2\ar[r]\ar[u]_{\rm exp}&\tder_2\ar[r]\ar[u]_{\rm exp}&\frak{un}(V)\ar[u]_{\rm exp}}\eqno(4.5)$$
where the upper left-hand group, $KRV_2$, is the prounipotent group actually
defined as the image in $SAut_2$ of 
$\kv\subset \sder_2$ under the exponential map, although the authors then go on to also provide a direct description 
of $KRV_2$ [AT, \S 5.1].
\vskip .3cm
Let us now recall the definition of the prounipotent group version $DS$ of the double shuffle Lie algebra 
$\ds$ originally given by Racinet in [R, Chap. 4, \S 1]; this is the group that Racinet denotes 
$DM_0({\bf k})$, but we take the base field ${\bf k}=\Q$; note that he also writes $\dm_0({\bf k})$ for $\ds$.

For any monomials $u,v\in \Q\langle x,y\rangle$, let the {\it shuffle product} 
$sh(u,v)\in \Q\langle x,y\rangle$ be defined recursively by
$$sh(1,u)=sh(u,1)=u, \ \ \ \ sh(Xu,Yv)=x\,sh(u,Yv)+y\,sh(Xu,v).\eqno(4.6)$$
It is well-known that the condition for a polynomial $f\in \Q\langle x,y\rangle$ to be a Lie polynomial
is equivalent to the condition 
$$\bigl(f|sh(u,v)\bigr)=0\eqno(4.7)$$ 
for all pairs of words $(u,v)$.  The elements
of the double shuffle Lie algebra $\ds$ are thus defined by (4.7) and the stuffle condition
$$\bigl(f_*|st(u,v)\bigr)=0\eqno(4.8)$$ 
for all words $u,v$ ending in $y$, where $f_*=\pi_y(f)+\sum_{n\ge 1}{{(-1)^{n-1}}\over{n}}(f|x^{n-1}y)y^n$
(cf. footnote to \S 1).

Let $DS$ be the group consisting of power series in $\Phi\in \Q\langle\langle x,y\rangle\rangle$ having 
constant term 1, no degree 1 or 2 terms, and satisfying two properties, which are
essentially group-like analogs of (4.7) and (4.8), namely
$$\bigl(\Phi|sh(u,v)\bigr)=\Phi(u)\Phi(v)\eqno(4.9)$$ 
for all pairs of words $(u,v)$ and
$$\bigl(\Phi_*|st(u,v)\bigr)=\Phi_*(u)\Phi_*(v)\eqno(4.10)$$ 
for all pairs of words $(u,v)$ both ending in $y$, where 
$$\Phi_*={\rm exp}\Bigl(\sum_{n\ge 1} {{(-1)^{n-1}}\over{n}}(\Phi|x^{n-1}y)y^n\Bigr)\pi_y(\Phi).$$

The elements of $\ds$ are Lie polynomials; as we saw in \S 1, the main result 
of [R] states that $\ds$ is a Lie algebra under the Poisson bracket (1.1). If $f\in \ds$,
then for any $g$ in the universal enveloping algebra ${\cal U}\ds$, the
multiplication in ${\cal U}\ds$ is given by the explicit formula $f\odot g
=fg+D_f(g)$.  Thus for $f\in \ds$, one can define $f^{\odot n}=f\odot
f^{\odot n-1}$, which gives an explicit polynomial formula for $f^{\odot n}$.
The exponential map of the Lie algebra is then given by
${\rm exp}^\odot(f)=\sum_{n\ge 0} {{1}\over{n!}}
f^{\odot n}$ as in (4.2).
In [R, Chap. 4, \S 3.3, corollaire 3.11] Racinet showed, using a method based on induction on the degree, that
$${\rm exp}^\odot(\ds)\simeq DS.\eqno(4.11)$$

The next theorem shows that there exists an injective map $DS\rightarrow KRV_2$, the group analog of the Lie 
algebra map of theorem 1.1. Given the results above on $KRV_2$ and $DS$, this is in fact nothing more than 
an immediate corollary of theorem 1.1.
\vskip .3cm
\noindent {\bf Theorem 4.2.} {\it There is an injective homomorphism of prounipotent groups
$DS\hookrightarrow KRV_2$ making the
following diagram commute:
$$\xymatrix{DS\ar@{^{(}->}[r]&KRV_2\\
\ds\ar@{^{(}->}[r]\ar[u]^{{\rm exp}^\odot}&\kv\ar[u]_{\rm exp}}$$}
\vskip .3cm
\noindent {\bf Proof.} 
Let $\rho:\ds\hookrightarrow \kv\subset \sder_2$ be the map of theorem 1.1; then 
the following diagram commutes by definition:
$$\xymatrix{{\rm exp}(\rho(\ds))\ar@{^{(}->}[r]&KRV_2\ar@{^{(}->}[r]&SAut_2\\
\rho(\ds)\ar[u]^{\rm exp}\ar@{^{(}->}[r]&\kv\ar[u]_{\rm exp}\ar@{^{(}->}[r]&\sder_2\ar[u]_{\rm exp},}\eqno(4.12)$$
where all the horizontal injections are just inclusions.

We also have the commutative diagram
$$\xymatrix{DS\ar[r]&{\rm exp}\bigl(\rho(\ds)\bigr)\\
\ds\ar[u]^{{\rm exp}^\odot}\ar[r]&\rho(\ds),\ar[u]_{\rm exp}}$$
where the left vertical arrow is Racinet's isomorphism (4.11), the right 
vertical arrow is the exponential isomorphism from (4.12), the bottom arrow is 
the isomorphism $\rho$ from theorem 1.1, and the top arrow is simply the
isomorphism defined by these other three arrows.  
Then the composition
$$DS\buildrel\sim\over\rightarrow {\rm exp}\bigl(\rho(\ds)\bigr)\subset KRV_2$$ 
is the desired injection.  \hfill{$\diamondsuit$}
\vskip .6cm
\noindent {\bf Acknowledgments.} Pierre Lochak and Samuel Baumard both
provided arguments for the second half of proposition 2.2,  the latter
being eventually used as it was shorter.  Much of the
spirit of the approach introduced here emerges from the reading of the
works of Jean Ecalle, who always insists that the situation must be
studied entirely via the symmetries that occur.  The terminology {\it anti}, 
{\it push} etc. is introduced purposely here with a view to eventually 
providing a more general introduction to his papers.  Finally, we warmly thank
the referee for a very detailed job with a great many useful suggestions and
corrections, in particular the addition of the final section of this paper.

\vfill\eject
\vskip 1cm
\noindent {\bf Appendix: Ecalle's theorem.}
\vskip .3cm
For $f\in \Q\langle x,y\rangle$, we write $f^r$ for the depth $r$ part of $f$ 
(i.e. the monomials containing exactly $r$ $y$'s), and 
$$f=\sum_{\e=(e_0,\ldots,e_r), r\ge 1} a_\e x^{e_0}y\cdots yx^{e_r}.\eqno(A.1)$$
To each $f\in \Q\langle x,y\rangle$, we associate two families of polynomials,
indexed by $\le r\le n$.  The first family, $vimo_f$, is a set of polynomials
in commutative variables $z_i$, and the second family, $ma_f$, is in commutative variables $u_i$.
$$vimo_f^r(z_0,\ldots,z_r)=\sum_{\e=(e_0,\ldots,e_r)} a_\e z_0^{e_0}\cdots z_r^{e_r},\eqno(A.2)$$
$$ma_f^r(u_1,\ldots,u_r)=vimo_f^r(0,u_1,u_1+u_2,\ldots,u_1+\cdots+u_r),\eqno(A.3)$$

Ecalle calls a {\it mould} any family $ma$ of functions $ma^r(u_1,\ldots,u_r)$, $r\ge 0$, with
$ma^0$ being a constant in a specified field.  He considers arbitrary functions, but 
in this appendix it is enough to consider only polynomial-valued moulds
$ma^r(u_1,\ldots,u_r)\in \Q[u_1,\ldots,u_r]$, with $ma^0=0$.  For any fixed integer $n\ge 1$,
such a mould is said to be {\it homogeneous of degree $n$} if $ma^r(u_1,\ldots,u_r)$ is a 
homogeneous polynomial of degree $n-r$.

Ecalle defines the following transformations of a mould $ma$ with $ma^0=0$:
$$swap(ma)^r(v_1,\ldots,v_r)=ma^r(v_r,v_{r-1}-v_r,\ldots,v_1-v_2)\eqno(A.4)$$
$$mantar(ma)^r(u_1,\ldots,u_r)=(-1)^{r-1}ma^r(u_r,\ldots,u_1)\eqno(A.5)$$
$$push(ma)^r(u_1,\ldots,u_r)=ma^r(-u_1-\cdots-u_r,u_1,\ldots,u_{r-1})\eqno(A.6)$$
$$teru(ma)^r(u_1,\ldots,u_r)=ma^r(u_1,\ldots,u_r)+{{1}\over{u_r}}
\bigl(ma^{r-1}(u_1,\ldots,u_{r-2},u_{r-1}+u_r)-ma^{r-1}(u_1,\ldots,u_{r-2},
u_{r-1})\bigr).\eqno(A.7)$$

The result of Ecalle that we use here is the following.
\vskip .3cm
\noindent {\bf Theorem A.1.} (Ecalle [E, \S 3.5, (3.64)]) {\it Let $n\ge 3$ and let
$\tilde f\in \ds_n$, so that $\tilde f$ is a homogeneous polynomial of degree $n$;
in particular $\tilde f^r=0$ if $r=0$ or $r\ge n$.  Let $ma$ be the mould $ma_{\tilde f}$ 
associated to $\tilde f$ as in (A.3).  Then $ma$ is a homogeneous mould of degree 
$n$, and for $1\le r\le n$, we have
$$teru(ma)^r=push\circ mantar\circ teru\circ mantar(ma)^r.
\eqno(A.8)$$}

The purpose of this appendix is to show that this theorem is equivalent to
theorem 3.3, by translating Ecalle's language back into terms of the
non-commutative variables $x$, $y$.  The first observation is that $mantar(ma)=ma$,
because $ma$ comes from a Lie polynomial. 
\vskip .3cm
\noindent {\bf Lemma A.2.} {\it Let $f\in {\rm Lie}_n[x,y]$ be a polynomial of homogeneous depth $r\ge 1$,
and let $ma$ be the mould associated to $f$ as in (A.3).  Then $mantar(ma)=ma$.}
\vskip .2cm
\noindent {\bf Proof.} Let $f$ be a polynomial of homogeneous degree $n\ge 3$ all of whose
terms of of fixed depth $r$; we write it as in (A.1) (with only the fixed value of $r$ giving non-zero
terms).  By the Lazard elimination theorem, any Lie polynomial
belongs to the polynomial ring generated by the polynomials $ad(x)^{i-1}(y)$ for $i\ge 1$.  Thus,
we can write
$$f=\sum_{{\bf c}} b_{{\bf c}}\ ad(x)^{c_1}(y)\cdots ad(x)^{c_r}(y).$$
We can show that we then have
$$ma_f^r(u_1,\ldots,u_r)=\sum_{{\bf c}} b_{{\bf c}}\ u_1^{c_1}\cdots u_r^{c_r};$$
in other words, the meaning of the coefficients of the mould $ma_f$ is that they reflect the
expression of $f$ as a polynomial in the $C_i$.  This idea was expressed by Racinet in [R] 
(Appendix A), but the proof is not given there.  It can be done by induction; the complete
proof is given in chapter 3 of the unpublished manuscript [S].

Now, if $P=ad(x)^{c-1}(y)$ and $anti(P)$ is as usual the polynomial 
obtained from $P$ by writing all its words backwards, then $anti(P)=(-1)^{c-1}P$. 
It follows that if $P$ is a product $P=ad(x)^{c_1-1}(y)\cdots ad(x)^{c_r-1}(y)$ and 
$P'=ad(x)^{c_r-1}(y)\cdots ad(x)^{c_1-1}(y)$, we have $P'=(-1)^{c_1+\cdots+c_r-r}anti(P)$.
Now assume that $f\in {\rm Lie}_n[x,y]$, so $(-1)^{n-1}anti(f)=f$.  This means that
$$\eqalign{f&=(-1)^{n-1}anti(f)=(-1)^{n-1}\sum_{{\bf c}} b_{{\bf c}}\ (-1)^{c_1+\cdots+c_r-r}ad(x)^{c_r-1}(y)\cdots ad(x)^{c_1-1}(y)\cr
&=(-1)^{r-1}\sum_{{\bf c}} b_{{\bf c}'}\ ad(x)^{c_1-1}(y)\cdots ad(x)^{c_r-1}(y),}$$
where if ${\bf c}=(c_1,\ldots,c_r)$, we write ${\bf c}'=(c_r,\ldots,c_1)$, so that
in particular $b_{{\bf c}'}=(-1)^{r-1}b_{{\bf c}}$.
Then
$$mantar(ma_f)^r(u_1,\ldots,u_r)=(-1)^{r-1}ma_f^r(u_r,\ldots,u_1)
=(-1)^{r-1}\sum_{{\bf c}} b_{{\bf c}}u_1^{c_r}\cdots u_r^{c_1}$$
$$=\sum_{{\bf c}} b_{{\bf c}'}u_1^{c_r}\cdots u_r^{c_1}
=\sum_{{\bf c}} b_{{\bf c}}u_1^{c_1}\cdots u_r^{c_r} =ma^r_f(u_1,\ldots,u_r).  $$
This concludes the proof.\hfill{$\diamondsuit$}
\vskip .3cm
The statement of Ecalle's theorem (A.8) for $r=1$ is easy to prove, since by (A.7), $teru(ma)^1(u_1)=ma^1(u_1)$,
and $push(ma^1(u_1))=ma^1(-u_1)$.  Now, if $n$ is even, it is well-known that if $\tilde f\in \ds_n$, then 
$\tilde f^1=0$, so $ma^1(u_1)=0$ and (A.8) holds.  If $n$ is odd, then either $\tilde f\in \ds_n$ also satisfies
$\tilde f^1=0$, so that again (A.8) holds, or $\tilde f^1=a\,ad(x)^{n-1}y$, in which case $ma^1(u_1)=au_1^{n-1}$,
so $push(ma^1(u_1))=ma^1(-u_1)=ma^1(u_1)$.  

Let us now give a reformulation of (A.8) for $2\le r\le n$.  By Lemma A.2, we can rewrite (A.8) as
$$swap\circ teru(ma)^r=swap\circ push\circ mantar\circ teru(ma)^r.\eqno(A.9)$$
The swap is obviously not necessary in the equality, but useful for the computation below
as it is easier to compute both sides as polynomials in the commutative variables $v_i$. 

By applying (A.4) to (A.7), we see that for $2\le r\le n$, the left-hand side is given by
$$swap\bigl(teru(ma)^r\bigr)(v_1,\ldots,v_r)=ma^r(v_r,v_{r-1}-v_r,\ldots,v_1-v_2)+$$
$${{1}\over{v_1-v_2}}\Bigl(ma^{r-1}(v_r,v_{r-1}-v_r,\ldots,v_3-v_4,v_1-v_3)-
ma^{r-1}(v_r,v_{r-1}-v_r,\ldots,v_3-v_4,v_2-v_3)\Bigr)$$
$$=vimo^r(0,v_r,\ldots,v_1)+{{1}\over{v_1-v_2}}\Bigl(vimo^{r-1}(0,v_r,\ldots,v_3,v_1)
-vimo^{r-1}(0,v_r,\ldots,v_3,v_2)\Bigr),\eqno(A.10)$$
where $vimo$ is the mould associated to $\tilde f$ as in (A.2).

Let us calculate the right-hand side of (A.9) one step at a time using (A.4)-(A.7).
$$\eqalign{swap&\circ push\circ mantar\circ teru(ma)^r=
swap\circ push\circ mantar\Bigl(ma^r(u_1,\ldots,u_r)\Bigr)\cr
&\quad +swap\circ push\circ mantar\Bigl({{1}\over{u_r}}\bigl(
ma^{r-1}(u_1,\ldots,u_{r-2},u_{r-1}+u_r)-ma^{r-1}(u_1,\ldots,u_{r-2},u_{r-1})\bigr)\Bigr)\cr
&=(-1)^{r-1}swap\circ push\Bigl(ma^r(u_r,\ldots,u_1)\Bigr)\cr
&\quad +(-1)^{r-1}swap\circ push\Bigl({{1}\over{u_1}}\bigl(
ma^{r-1}(u_r,\ldots,u_{3},u_{1}+u_2)-ma^{r-1}(u_r,\ldots,u_3,u_2)\bigr)\Bigr)\cr
&=(-1)^{r-1}swap\Bigl(ma^r(u_{r-1},\ldots,u_2,u_1,-u_1-\cdots-u_r)\Bigr)\cr
&\quad +(-1)^{r-1}swap\Bigl({{1}\over{(-u_1-\cdots-u_r)}}\bigl(
ma^{r-1}(u_{r-1},\ldots,u_{2},-u_2-\cdots-u_r)-ma^{r-1}(u_{r-1},\ldots,u_2,u_1)\bigr)\Bigr)\cr
&=(-1)^{r-1}ma^r(v_2-v_3,\ldots,v_{r-1}-v_r,v_r,-v_1)\cr
&\quad +(-1)^{r-1}{{1}\over{-v_1}}\bigl(ma^{r-1}(v_2-v_3,\ldots,v_{r-1}-v_r,v_r-v_1)
-ma^{r-1}(v_2-v_3,\ldots,v_{r-1}-v_r,v_r)\bigr)\cr
&=(-1)^{r-1}vimo^r(0,v_2-v_3,\ldots,v_2-v_r,v_2,v_2-v_1)\cr
&\quad +{{(-1)^r}\over{v_1}}\Bigl(vimo^{r-1}(0,v_2-v_3,\ldots,v_2-v_r,v_2-v_1)-
vimo^{r-1}(0,v_2-v_3,\ldots,v_2-v_r,v_2)\Bigr).}\eqno(A.11)$$
The following useful elementary identities will simplify the form of (A.11): for any 
$vimo$ associated to a polynomial as in (A.2), we have
$$vimo^r(z_0,\ldots,z_r)=(-1)^{n-r}vimo^r(-z_0,\ldots,-z_r),\eqno(A.12)$$
and if $vimo$ is associated to a Lie polynomial, then
$$vimo^r(z_0,z_1,\ldots,z_r)=vimo^r(0,z_1-z_0,\ldots,z_r-z_0).\eqno(A.13)$$
Note that the meaning of (A.13) is that any value (called $z_0$) can be added to each
argument of $vimo^r$ without changing the value of the function.  Let us quickly indicate
the easy proof of (A.13) by induction.  For $r=1$, up to scalar multiple, we must have 
$$f^1=ad(x)^m(y)=\sum_{i=0}^m (-1)^i\bigl({{m}\atop{i}}\bigr)x^{m-i}yx^i,\ \ \ {\rm so}\ \ \ 
vimo^1(z_0,z_1)=\sum_{i=0}^m (-1)^i\bigl({{m}\atop{i}}\bigr)z_0^{m-i}z_1^i,$$
which is equal to $(-1)^m(z_1-z_0)^m=vimo^1(0,z_1-z_0)$.  Now assume that (A.13) holds
up to depth $r-1$ and consider a Lie polynomial $f$ of homogeneous depth $r$.  By linearity,
we may assume that $f=[g,h]$, where $g$ and $h$ are homogeneous depths $s<r$ and $t<r$ respectively,
with $r=s+t$.  Then we have
$$vimo_f^r(z_0,\ldots,z_r)=vimo_g^s(z_0,\ldots,z_s)vimo_h^t(z_s,\ldots,z_{s+t})-
vimo_h^t(z_0,\ldots,z_t)vimo_g^s(z_t,\ldots,z_{s+t}),$$
so using repeated applications of (A.13) to the $vimo_g$ and $vimo_h$ factors by the induction hypothesis,
we have 
$$\eqalign{vimo_f^r(&0,z_1-z_0,\ldots,z_r-z_0)
=vimo_g^s(0,z_1-z_0,\ldots,z_s-z_0)vimo_h^t(z_s-z_0, \ldots,z_{s+t}-z_0)-\cr
&\qquad\qquad \qquad\qquad \qquad\qquad \qquad\qquad 
vimo_h^t(0,z_1-z_0,\ldots,z_t-z_0)vimo_g^s(z_t-z_0,\ldots,z_{s+t}-z_0)\cr
&=vimo_g^s(z_0,z_1,\ldots,z_s)vimo_h^t(0,z_{s+1}-z_s,\ldots,z_{s+t}-z_s)- \cr
&\qquad\qquad \qquad\qquad \qquad\qquad \qquad\qquad 
vimo_h^t(z_0,z_1,\ldots,z_t)vimo_g^s(0,z_{t+1}-z_t,\ldots,z_{s+t}-z_t)\cr
&=vimo_g^s(z_0,z_1,\ldots,z_s)vimo_h^t(z_s,z_{s+1},\ldots,z_{s+t})- 
vimo_h^t(z_0,z_1,\ldots,z_t)vimo_g^s(z_t,z_{t+1},\ldots,z_{s+t})\cr
&=vimo_f^r(z_0,\ldots,z_r).}$$
This yields (A.13).
\vskip .3cm
Now, applying (A.12) to (A.11) yields
$$(-1)^{n-1}vimo^r(0,v_3-v_2,\ldots,v_r-v_2,-v_2,v_1-v_2)$$
$$+{{(-1)^n}\over{v_1}}\Bigl(vimo^{r-1}(0,v_3-v_2,\ldots,v_r-v_2,v_1-v_2)-
vimo^{r-1}(0,v_3-v_2,\ldots,v_r-v_2,-v_2)\Bigr).$$
Applying (A.13) to this with $z_0=v_2$, i.e. adding $v_2$ to each argument of $vimo^r$, then yields
$$ (-1)^{n-1}\Bigl[vimo^r(v_2,v_3,\ldots,v_r,0,v_1)-
{{1}\over{v_1}}\Bigl(vimo^{r-1}(v_2,v_3,\ldots,v_r,v_1)-vimo^{r-1}(v_2,v_3,\ldots,
v_r,0)\Bigr)\Bigr].\eqno(A.14)$$
Since if $vimo$ is the mould associated to a polynomial $\tilde f\in \ds_n$, then
$vimo^0=vimo^n=0$, Ecalle's theorem can be expressed by the equalities (A.10)=(A.14) 
for $1\le r\le n$.
\vskip .2cm
Let us now show that the statement of theorem 3.3 can be deduced from the equalities
(A.10)=(A.14) for $2\le r\le n$.
\vskip .3cm
\noindent {\bf Proposition A.3.} {\it Let $\tilde f\in\ds_n$ for $n\ge 3$, and write
$\tilde f=\tilde f_xx
+\tilde f_yy$.  Then $\tilde f_x+\tilde f_y$ is antipalindromic.}
\vskip .2cm
\noindent {\bf Proof.}  We will show the identity
$$\bigl(\tilde f_x+\tilde f_y\bigr)^r=(-1)^{n-1}anti\bigl(\tilde f_x+\tilde f_y\bigr)^r,$$ 
separately for each depth $1\le r\le n-1$ occurring in $\tilde f$.

These equalities are equivalent to the equalities of polynomials in commutative variables
$$vimo_{\tilde f_x^r+\tilde f_y^r}(z_0,\ldots,z_r)=(-1)^{n-1}vimo^r_{anti(\tilde f_x^r+
\tilde f_y^r)}(z_0,\ldots,z_r)\eqno(A.16)$$
for $1\le r\le n-1$.

To prove the proposition, we will deduce (A.16) from 
Ecalle's theorem, i.e. from the set of equalities (A.10)=(A.14) for $2\le r\le n$.
To do this, we explicitly compute both sides of (A.16).

Each term of the polynomial $(\tilde f_x+\tilde f_y)^r$ comes either from a term in 
$\tilde f^r$ ending with $x$ (i.e. from $\tilde f_x^rx$) or from a term in $\tilde f^{r+1}$ ending with $y$
(i.e. from $\tilde f_y^ry$), by cutting off the final letter.  Let us first find the 
$vimo$ polynomials associated to $\tilde f_y^ry$ and $\tilde f_x^rx$.  

Let $\tilde f^{r+1}=\sum_{\e=(e_0,\ldots,e_{r+1})} a_\e\, x^{e_0}y\cdots yx^{e_{r+1}}.$
Since $\tilde f^{r+1}$ is homogeneous in depth $r+1$, we have
$$vimo_{\tilde f^{r+1}}(z_0,\ldots,z_{r+1})=
\sum_{\e=(e_0,\ldots,e_{r+1})} a_\e z_0^{e_0}\cdots z_{r+1}^{e_{r+1}}.$$
Similarly, writing 
$\tilde f_y^ry=\sum_{\e=(e_0,\ldots,e_r,0)} a_\e\,x^{e_0}y\cdots x^{e_r}y$,
the depth $r+1$ polynomial $vimo_{\tilde f_y^ry}$ is given by
$$vimo_{\tilde f_y^ry}(z_0,\ldots,z_{r+1})=
\sum_{\e=(e_0,\ldots,e_r,0)} a_\e z_0^{e_0}\cdots z_r^{e_r}=
vimo_{\tilde f^{r+1}}(z_0,\ldots,z_r,0).$$
Since we have $\tilde f_y^r=\sum_{\e=(e_0,\ldots,e_r)} a_\e\,x^{e_0}y\cdots yx^{e_r}$, we see
that $vimo_{\tilde f_y^r}(z_0,\ldots,z_r)=vimo_{\tilde f_y^ry}(z_0,\ldots,z_{r+1})$, i.e.
$$vimo_{\tilde f_y^r}(z_0,\ldots,z_r)=vimo_{\tilde f^r_yy}(z_0,\ldots,z_{r+1})=vimo_{\tilde f^{r+1}}(z_0,
\ldots,z_r,0).\eqno(A.17)$$
To find the $vimo$ associated to $\tilde f_x^rx$, we consider this polynomial
as the difference $\tilde f_x^rx=\tilde f^r-\tilde f^{r-1}_yy$.  Thus, using (A.17) for $r-1$ instead
of $r$, we have
$$vimo_{\tilde f^r_xx}(z_0,\ldots,z_r)=vimo_{\tilde f^r-\tilde f^{r-1}_yy}(z_0,\ldots,z_r)=
vimo_{\tilde f^r}(z_0,\ldots,z_{r-1},z_r)-vimo_{\tilde f^r}(z_0,\ldots,z_{r-1},0).\eqno(A.18)$$
Because we know that there is an $x$ at the end of every word of the polynomial
$f^r_xx$, this polynomial is divisible by $z_r$, and we have
$$vimo_{\tilde f^r_x}(z_0,\ldots,z_r)=
{{1}\over{z_r}}\Bigl(vimo_{\tilde f^r}(z_0,\ldots,z_{r-1},z_r)-vimo_{\tilde f^r}(z_0,\ldots,z_{r-1},0)
\Bigr).\eqno(A.19)$$
Putting (A.17) and (A.19) together yields the following expression for the 
left-hand side of the desired equality (A.16):
$$vimo_{(\tilde f_x+\tilde f_y)^r}(z_0,\ldots,z_r)= vimo_{\tilde f^{r+1}}(z_0, \ldots,z_r,0)
+{{1}\over{z_r}}\Bigl(vimo_{\tilde f^r}(z_0,\ldots,z_{r-1},z_r)-vimo_{\tilde f^r}(z_0,\ldots,z_{r-1},0)
\Bigr).\eqno(A.20)$$
Since $anti$ corresponds to reversing the order of $z_0,\ldots,z_r$, the right-hand side of (A.16) is then 
given by
$$(-1)^{n-1}vimo_{anti((\tilde f_x+\tilde f_y)^r)}(z_0,\ldots,z_r)= $$
$$(-1)^{n-1}
\Bigl[vimo_{\tilde f^{r+1}}(z_r, \ldots,z_0,0)
+{{1}\over{z_0}}\Bigl(vimo_{\tilde f^r}(z_r,\ldots,z_1,z_0)-vimo_{\tilde f^r}(z_r,\ldots,z_1,0)
\Bigr)\Bigr],\eqno(A.21)$$
so the statement of the proposition is equivalent to the set of equalities
(A.20)=(A.21) for $1\le r\le n-1$.
\vskip .2cm
Thus it remains only to show that Ecalle's set of equalities (A.10)=(A.14) for $2\le r\le n$
implies the set of equalities (A.20)=(A.21) for $1\le r\le n-1$.  By (A.13), we can add the same
quantity to every argument of $vimo$ and not change its value, so we first use this to rewrite
(A.10), by adding the uantity $-v_1$ to every argument of the three $vimo$ terms in (A.10):
$$(A.10)=vimo_{\tilde f^r}(-v_1,v_r-v_1,\ldots,v_2-v_1,0)$$
$$+{{1}\over{v_1-v_2}}
\Bigl(vimo_{\tilde f^{r-1}}(-v_1,v_r-v_1,\ldots,v_3-v_1,0)
-vimo_{\tilde f^{r-1}}(-v_1,v_r-v_1,\ldots,v_3-v_1,v_2-v_1)\Bigr).$$
Now we apply the variable change 
$$z_0=-v_1,\ z_1=v_r-v_1,\ldots,z_{r-1}=v_2-v_1\eqno(A.22)$$
to this, to obtain
$$=vimo_{f^r}(z_0,z_1,\ldots,z_{r-1},0)+{{1}\over{-z_{r-1}}}
\Bigl(vimo_{f^{r-1}}(z_0,z_1,\ldots,z_{r-2},0)
-vimo_{f^{r-1}}(z_0,\ldots,z_{r-2},z_{r-1})\Bigr).\eqno(A.23)$$
This is equivalent to (A.20), for $r-1$ instead of $r$.

Next, we use (A.13) to rewrite (A.14), adding the quantity $-v_1$ to every argument in the three
$vimo$ terms that appear in (A.14):
$$ (-1)^{n-1}\Bigl[vimo_{\tilde f^r}(v_2-v_1,v_3-v_1,\ldots,v_r-v_1,-v_1,0)$$
$$- {{1}\over{v_1}}\Bigl(vimo_{\tilde f^{r-1}}(v_2-v_1,v_3-v_1,\ldots,v_r-v_1,0)
-vimo_{\tilde f^{r-1}}(v_2-v_1,v_3-v_1,\ldots, v_r-v_1,-v_1)\Bigr)\Bigr],$$
and then the variable change (A.22), which yields
$$ (-1)^{n-1}\Bigl[vimo_{\tilde f^r}(z_{r-1},z_{r-2},\ldots,z_1,z_0,0)$$
$$- {{1}\over{z_0}}\Bigl(vimo_{\tilde f^{r-1}}(z_{r-1},z_{r-2},\ldots,z_1,0)
-vimo_{\tilde f^{r-1}}(z_{r-1},z_{r-2},\ldots, z_1,z_0)\Bigr)\Bigr].$$
This is exactly (A.21) for $r-1$ instead of $r$.  Thus Ecalle's equalities 
(A.10)=(A.14) for $2\le r\le n$ imply the desired equalities (A.16) for $1\le r\le n-1$
as desired.\hfill{$\diamondsuit$}
\vfill\eject

\vfill\eject
\noindent {\bf References}
\vskip .3cm
\noindent
[AT] A. Alekseev, C. Torossian, The Kashiwara-Vergne Conjecture and Drinfeld's 
associators, arXiv:0802.4300, preprint 2009.
\vskip .3cm
\noindent
[B] N. Bourbaki, {\it Groupes et alg\`ebres de Lie, Chapitres 2 et 3}, Springer-Verlag, Berlin
1982.
\vskip .3cm
\noindent
[CS] S. Carr, L. Schneps, Combinatorics of the double shuffle Lie algebra,
to appear in {\it Grothendieck-Teichm\"uller theory and arithmetic
geometry}, Proceedings of the 3rd MSJ-SI conference in Kyoto, October 2010.
\vskip .3cm
\noindent
[E] J. Ecalle, The flexion structure and dimorphy: flexion units, singulators, generators,
and the enumeration of multizeta irreducibles, to appear in Ann. Scuo. Norm. Pisa, 2011.
\vskip .3cm
\noindent
[F] H. Furusho, Double shuffle relation for associators, to appear in {\it 
Annals of Math.}, 2011.
\vskip .3cm
\noindent
[IKZ] K. Ihara, M. Kaneko, D. Zagier, Derivation and double shuffle relations
for multiple zeta values, {\it Compos. Math.} {\bf 142} (2006), no. 2,
307-338. 
\vskip .3cm
\noindent
[R] G. Racinet, S\'eries g\'en\'eratrices non-commutatives de polyz\^etas
et associateurs de Drinfel'd, Ph.D. thesis, 2000.
\vskip .3cm
\noindent
[S] L. Schneps, ARI, GARI, Zig and Zag: Ecalle's theory of multiple zeta values, book manuscript, 2011.
\vskip 2cm
\noindent 
Tel: (33) 1 44 27 53 55
\vskip .1cm
\noindent
Email: leila@math.jussieu.fr
\vskip .1cm
\noindent 
Institut de Math\'ematiques de Jussieu, 4 place Jussieu, Case 247, 75252 Paris Cedex, France
\vskip .2cm

\bye